\documentclass[11pt]{article}
\usepackage{amssymb,amsmath,epsfig}

\newtheorem{theorem}{Theorem}[section]
\newtheorem{thm}[theorem]{Theorem}
\newtheorem{Thm}{Theorem}

\newtheorem{prop}[theorem]{Proposition}

\newtheorem{lemma}[theorem]{Lemma}
\newtheorem{cor}[theorem]{Corollary}

\newenvironment{proof}{\par\medskip\noindent{\em Proof. }}{\hfill $\square$\par\medskip}

\def\<{\langle}
\def\>{\rangle}
\def\Z{\mathbb{Z}}

\def\Q{\mathbb{Q}}
\def\N{\mathbb{N}}
\def\T{\mathbb{T}}
\def\G{\Gamma}

\def\o{\omega}

\title{Normalizers in Limit Groups}
\author{Martin R. Bridson and James Howie}
\begin{document}

\maketitle
\begin{abstract} Let $\G$ be a limit group,
 $S\subset\G$ a subgroup, and $N$ the normaliser
 of $S$. If $H_1(S,\mathbb Q)$ has
finite $\Q$-dimension, then $S$ is finitely generated and
either $N/S$ is finite or $N$ is abelian.
This result has applications to the study of subdirect
products of limit groups.
\end{abstract}

Limit groups (or {\em $\o$-residually free groups})
have received a good deal of attention in recent
 years, primarily due to the groundbreaking
work of Z.~Sela (\cite{S1,S2} {\em{et seq.}}).
See for example \cite{AB,BF,CG,Guirardel,Paulin}.

O.~Kharlampovich, A.~Myasnikov \cite{KM1,KM2}
and others (see, for example, \cite{FGMRS,Kapovich,KMRS}) have studied limit
groups extensively
under the more traditional name of 
{\em{fully residually free groups}}, which
appears to have been first 
introduced by B.~Baumslag in \cite{Ba}. This name
reflects the fact that limit groups are precisely
those finitely
generated groups $\G$ such that
 for each finite subset $T\subset \G$ there exists a 
 homomorphism
from $L$ to a free group that is injective on $T$.

Examples of limit groups include all finitely generated free 
or free abelian groups, and the fundamental groups
of all closed surfaces of Euler characteristic at most $-2$.
The free product of finitely many limit groups is again a limit group,
which leads to further examples. More sophisticated examples
are described in some of the articles cited above.

The purpose of this note is to  contribute some
results on the subgroup structure of limit groups. 

\begin{Thm}\label{norm} If $\G$ is a limit group
and $H\subset\G$ is a finitely generated subgroup,
then either $H$ has finite index in its normaliser
or else the normaliser of $H$ is abelian.
\end{Thm}

This result is in keeping with the expectation that
finitely generated subgroups of limit groups should
be quasi-isometrically embedded.

We shall use the following result to circumvent the
difficulty that {\em{a priori}} one does not know that
normalisers in limit groups are finitely generated.

\begin{Thm}\label{1factor}
 Let $\G$ be a limit group and  $S\subset \G$  a subgroup. 
 If $H_1(S,\Q)$ has finite $\Q$-dimension,
then  $S$ is finitely generated (and hence is a limit group).
\end{Thm}

Theorems \ref{norm} and \ref{1factor} combine to give the result stated
in the abstract.
Theorem \ref{norm} plays an important
role in our work on the subdirect products of limit
groups \cite{BH,BH2}. In the present note, we content ourselves
with noting the following easy consequence of Theorem \ref{norm}.

\begin{Thm}\label{split}
Suppose that $\G_1,\dots,\G_n$ are limit groups 
and let $$S\subset \G_1\times\cdots\times \G_n$$ be an arbitrary subgroup.
If $L_i=\G_i\cap S$ is non-trivial and 
finitely generated for $i\le r$, then
a subgroup $S_0\subset S$ of finite index 
splits as a direct product
$$
S_0= L_1\times\dots\times L_r\times (S_0\cap G_r),
$$
where $G_r = \G_{r+1}\times\dots\times \G_n$.
\end{Thm}

Thus, up to finite index, the study of subdirect products 
of limit
groups reduces to the case where the normal subgroups 
$L_i\triangleleft S$ are all infinitely generated.

The remainder of this paper is structured as follows.  
In section
\ref{lgt} below we recall the definitions and record
 some elementary
properties of limit groups and $\o$-residually free towers. 
In section
\ref{nsb} we prove a result about normal subgroups 
of groups acting
on trees. In section \ref{rsn} we describe a large
class of groups with the property that each 
finitely generated, non-trivial, normal subgroup is
of finite index. This class includes all non-abelian
limit groups. In sections \ref{homfin}, 
\ref{nfg} and \ref{spl} we
prove Theorems \ref{1factor}, \ref{norm} and 
\ref{split} respectively.

\section{Limit groups and towers}\label{lgt}

Our results rely  on  the fact
 that limit groups
are the finitely generated
subgroups of $\omega$-residually free tower groups
(\cite{S2} and
\cite{KM2}).

\medskip\noindent{\bf Definition} 
An {\em $\omega$-rft space} of {\em{height}} $h\in\N$
is  defined by induction on $h$. An {\em $\omega$-rft {group}}
is the fundamental group of an $\omega$-rft space.

A height $0$ tower is the wedge
(1-point union) of a finite collection of circles, closed hyperbolic
surfaces and tori $\T^n$ (of arbitrary dimension), except that the closed
surface of Euler characteristic $-1$ is excluded.

An $\omega$-rft space $Y_n$ of height $h$ is obtained from 
an $\omega$-rft space $Y_{h-1}$ of  height $h-1$ by attaching a
{\em{block}} of one of the following types:

\noindent{\em{Abelian:}} $Y_h$ is obtained from $Y_{h-1}\sqcup\T^m$
by identifying a coordinate circle in $\T^m$ with any loop $c$ in 
$Y_{h-1}$ such that $\langle c\rangle\cong\Z$ is a maximal abelian subgroup
of $\pi_1Y_{h-1}$.

\noindent{\em{Quadratic:}} One takes a connected, compact surface $\Sigma$
 that is either a punctured torus or has Euler characteristic at most $-2$, and
obtains $Y_h$  from $Y_{h-1}\sqcup\Sigma$ by identifying each boundary
component of $\Sigma$ with a homotopically non-trivial loop in $Y_{h-1}$; these
identifications must be chosen so that there exists a  retraction 
$r:Y_h\to Y_{h-1}$, and $r_*(\pi_1\Sigma)\subset\pi_1Y_{h-1}$ must be non-abelian.

\medskip

\begin{thm}\label{embed}{\em\cite{S2}} A given group is a 
limit group
if and only if it is isomorphic to a
finitely generated subgroup of an $\o$-rft group.
\end{thm}

A useful sketch of the proof can be found in \cite[Appendix]{AB}.

This is a powerful theorem that allows one to prove many
interesting facts about limit groups by induction on  height.

The {\em{height}} 
of a limit group $S$ is the minimal height of an 
$\o$-rft group
that has a subgroup isomorphic to $S$. Limit
groups of height $0$ are free products of 
finitely many free abelian groups and of surface
groups of Euler characteristic at most $-2$.

 The splitting described in the following
 lemma is that which the Seifert-van Kampen
Theorem associates to the addition of the
final block in the tower construction.

\begin{lemma}\label{l:basic}
The fundamental group of an $\o$-rft  
 space of height $h\ge 1$ splits as a 2-vertex graph 
of groups: one
of the vertex groups $A_1$ is the fundamental
group of an $\o$-rft  
 space  of height $h-1$ and the other $A_2$
is a free or free-abelian group
of finite rank at least 2; the edge groups are
 maximal infinite cyclic subgroups
of $A_2$. If $A_2$ is abelian then there is only
one edge group and this is a maximal abelian subgroup
of $A_1$.
\end{lemma}

Since an arbitrary limit group (and hence an arbitrary subgroup
of a limit group) is a subgroup of
an $\o$-rft group, one can apply  Bass-Serre theory 
to deduce:

\begin{lemma}\label{unconflate}
If $S$ is a subgroup of a limit group of height $h\ge 1$, then
$S$ is  the fundamental group of a 
bipartite graph of groups in which the
edge groups are cyclic; the vertex groups fall into two types 
corresponding to the bipartite partition of vertices:
type (i) vertex groups are isomorphic to subgroups of 
a limit group of height $(h-1)$;
type (ii) vertex groups are all free or
all free-abelian.
\end{lemma}

Restricting attention to limit groups
(rather than arbitrary subgroups)
we will also use the following variant of the above decomposition.

\begin{lemma}\label{BSsub} 
If $\G$ is a freely indecomposable
limit group of height $h\ge 1$, then
it is  the fundamental group of a finite
graph of groups that has infinite
cyclic edge groups and 
has a 
vertex group that is a non-abelian
limit group of height $\le h-1$.
\end{lemma}

\begin{proof}
Every limit group arises as a finitely generated subgroup of an
$\o$-rft group of the same height.  
Suppose $L$ is an $\o$-rft group of height $h$ 
that contains $\G$ as a subgroup. Since $\G$ is finitely
generated, the graph of groups structure given 
by Lemma \ref{unconflate} 
can be replaced by the
finite core graph $\mathcal C$ obtained by taking 
the quotient by $\G$ of a
{\em minimal}
 $\G$-invariant subtree of the Bass-Serre tree for $L$.
The vertex groups of $\mathcal C$ are
finitely generated, and since $\G$ does not have height
$\le h-1$, there is at least one vertex group of each
type (in the
terminology of Lemma \ref{unconflate}). As $\G$
is freely indecomposable, the edge groups in $\mathcal C$
 are infinite cyclic.

If the
block added at the top of the tower for $L$ is quadratic,
then 
each vertex group $\G_v$ of type (ii) in $\mathcal C$
must be a non-abelian free group (and hence a non-abelian
limit group of height $0\le h-1$). Indeed,
 distinct conjugates of the edge groups incident at 
the corresponding vertex group in $L$ are maximal cyclic subgroups 
of the vertex group and generate a non-abelian 
free subgroup. Thus if $\G_v$ were cyclic, then there would
be a single incident edge $e$ at $v$ and $\G_e\to\G_v$ would be
an isomorphism, contradicting the fact that $\mathcal C$
was chosen to be  minimal.

To complete the proof we shall argue that
if the block added at the top of the 
tower for $L$ is abelian, then each vertex of type
(i) in $\mathcal C$ is non-abelian. 

We have
$L=A_1\ast_{\langle \zeta\rangle}A_2$, 
as in Lemma \ref{l:basic}, where 
$\langle \zeta\rangle\cong \mathbb Z$.
In $\mathcal C$,  the edge groups incident at each
vertex of type (i) (which are conjugates of subgroups of 
$\langle \zeta\rangle$) intersect trivially and 
are maximal abelian. 
Thus if a vertex group $\Gamma_v$ of type (i)
were abelian, then it would
be cyclic, there would be only one incident edge $e$, and
the inclusion $\G_e\to\G_v$ would be an isomorphism.
As above, this would contradict the
minimality of $\mathcal C$.
\end{proof}

\section{Normal Subgroups in Bass-Serre Theory}\label{nsb}

As was evident in the previous section, we are assuming that
 the reader is familiar with the basic theory of groups
acting on trees.  We  need the  following fact
from this theory.

\begin{prop}\label{BSp} If a finitely
generated group $\G$ acts by isometries on a tree $T$,
 then either $\G$ fixes a point of $T$
  or else $\G$ contains a
hyperbolic isometry.

In the latter case, the union of the axes of 
the hyperbolic elements
is the unique minimal $\G$-invariant subtree of $T$.
\end{prop}

The following dichotomy for normal subgroups in Bass-Serre theory
ought to be well known but we are unaware of any appearance of it
in the literature.
(The results of \cite{Moon} are similar in spirit but different in
content.)

The first of the stated outcomes arises, for example,
when $\G=G_0\ast_AG_1$ and  $A$ is normal
in both $G_0$ and $G_1$. The second outcome arises, for example,
when $\G$ is the amalgamation
of $N_0\rtimes A$ and $N_1\rtimes A$, and $N=N_0\ast N_1$.

\begin{prop}\label{Nfg} If the group $\G$ splits over the
subgroup $A\subset\G$ 
and the normal subgroup $N\triangleleft\G$ is finitely
generated, then either $N\subset A$  
or else $AN$ has finite index in $\G$.
\end{prop}

\begin{proof} 
Consider the action of $N$ on the Bass-Serre tree $T$
associated to the splitting of $\G$ over $A$.  
If the fixed-point set of  $N$ is non-empty, then it is the
whole of $T$, because it is $\G$-invariant and the action 
is minimal.
Thus $N$ lies in the edge-stabiliser $A$.  

If $N$ does not have a fixed point, then by
the preceding proposition it must
contain  hyperbolic isometries and
the union of the axes of its hyperbolic elements is 
the minimal $N$-invariant subtree of $T$.
Since $N$ is normal this tree is also $\G$-invariant, 
hence equal to $T$.

The edges of $T$ are
indexed by the cosets $\G/A$ and the action of $g\in\G$ is
$g.\gamma A = (g\gamma)A$. Thus the edges of the graph $X=N\backslash T$
are indexed by the double cosets $N\backslash\G/A$. 
The quotient $\G/N$ acts on
$X$ and the stabilizer of the edge indexed $NgA$ is the image in $\G/N$
of $gAg^{-1}$. It follows that some (hence every)
edge in $X$ has a finite orbit if and only
the image of $A$ in $\G/N$ is of finite index.

But $X$ is a finite graph. Indeed, since $N$ is finitely generated, 
there is a compact subgraph $Y\subset X$
such that the inclusion of the corresponding graph of groups induces
an isomorphism of fundamental groups, which implies that the action
of $N$ preserves the (connected) preimage of $Y$ in $T$.
Since the action of $N$ on $T$ is minimal, this preimage is the whole of $T$.
\end{proof} 

\section{A class of groups with restricted normal
subgroups}\label{rsn}

In this section we describe a large class 
of groups that have the property
that each of their non-trivial, finitely generated, 
normal subgroups is of finite index. This leads to
the following special case of Theorem \ref{norm}.

\begin{thm} \label{Lfi} Let $\G$ be a non-abelian 
limit group and let
$N\neq 1$ be a normal subgroup.
If $N$
is finitely generated, then $\G/N$ is finite.
\end{thm}

Let $\mathcal J$ denote the class of groups that have no proper infinite
quotients. (Thus $\mathcal J$ is the union of the class of finite groups
and the class of just-infinite groups.)
Let $\mathcal G$ be the smallest class of groups 
that contains $\Z$ and all non-trivial free products of groups,
and satisfies the
following condition:
\begin{enumerate}\item[$\ast$]
If $H$ is in $\mathcal G$  and $A$ is in $\mathcal J$, and if $A$ has infinite
index in $H$,  then
any amalgamated free product of the form $H\ast_A B$ is in
$\mathcal G$, and so is any HNN extension of the form $H\ast_A$.
\end{enumerate}
 
\begin{thm} \label{fi} 
If $\G$ is in $\mathcal G$, 
then every finitely generated, non-trivial, normal
subgroup  in $\G$ is of finite index. 
\end{thm}

\begin{proof} Given $\G$, we argue by induction on 
its {\em{level}} in  $\mathcal G$, namely
the number of steps in which $\G$ is constructed 
by taking amalgamated free products and HNN extensions as in
$(\ast)$. At level 0 we have $\G\cong\Z$ or $\G$ a non-trivial free product. 
In the former case the result is clear; in the latter it follows from
Proposition \ref{Nfg} with $A=\{1\}$.

Assume that $\G= H\ast_A B$ or $\G= H\ast_A$, where every finitely generated
normal subgroup of  $H$ has finite index and $A\in\mathcal J$ has infinite index
in $H$. Let $N\subset\G$ be finitely generated and normal. According
to Proposition \ref{Nfg}, either $N\subset A$ or else $NA$ has finite index in $\G$.
The former case is impossible because it would imply  $N\subset H$ and
 by induction $|H/N|$ would be finite, implying
that $|H/A|$ is finite, contrary to hypothesis. 
Thus the image of $A$ has finite index in $\G/N$. 
 Since $A$ has no proper infinite quotients,
it follows that either $N$ has finite index in $\G$ (as desired), or else
$A$ is infinite and $N\cap A=\{1\}$. 

We claim that this  last possibility is absurd. First note that since
 $A$ has infinite index in $H$ and projects to a subgroup of finite
index in $\G/N$, we know that  $H\cap N$ is non-trivial. Secondly,
because $A$ (and hence each of its conjugates
in $\G$) intersects $N$ trivially, the quotient of the Bass-Serre tree
for $\G$ gives a graph of groups decomposition for $N$ with trivial
edge groups. One of the vertex groups in this decomposition is $H\cap N$.
Thus $H\cap N$ is a free factor of the finitely generated group $N$, and
hence is finitely generated. By induction, if follows that $|H/(H\cap N)|$ is finite.
Hence $|A/(A\cap N)|$ is finite, which is the contradiction we 
seek.
\end{proof}

\begin{prop}\label{LG}
All non-abelian limit groups lie in $\mathcal G$.
\end{prop}

\begin{proof} This follows easily from the fact that an 
arbitrary limit group $\G$ is
a subgroup of an $\o$-rft group. A non-abelian limit
group of height $0$ is either a non-trivial free product or a
surface group of Euler characteristic at most $-2$ (which therefore
splits as $F\ast_\Z$ with $F$ free).  Such groups
clearly lie in $\mathcal G$,
so we may assume that the tower 
has height $h\ge 1$ and that $\G$ is freely indecomposable.

As in Lemma \ref{BSsub}, we  decompose $\G$ as
$\pi_1\mathcal C$. That
lemma singles out a non-abelian vertex group $\Gamma_v$
that is 
a limit group of height $\le h-1$. By induction $\Gamma_v$ lies in $\mathcal G$. By conflating the remaining vertices 
of $\mathcal C$ to a single vertex, we can express 
$\Gamma$ as the fundamental group of a 2-vertex
graph of groups where one vertex group is $\Gamma_v$
and the other is the fundamental group of the
full subgraph of groups spanned by the remaining
vertices of $\mathcal C$; the edge groups are cyclic.
A secondary induction on
the number of edges in the decomposition
completes the proof. 
\end{proof}

If $\G$ is a hyperbolic limit group then it has cohomological
dimension at most 2, so by a theorem of Bieri \cite{bieri2} its finitely
presented normal subgroups are free or of finite index. Since any
finitely generated
subgroup of a limit group is a limit group, and limit groups are
finitely presented,  the novel content of Theorem \ref{Lfi}
in the hyperbolic case is essentially contained in: 

\begin{cor} If $F$ is a finitely generated, non-abelian free group, then no
group of the form $F\rtimes\Z$ is a limit group.
\end{cor}

\section{Homological finiteness and finite generation}\label{homfin}

The main object of this section is to
 prove Theorem \ref{1factor}, which is
required in order to prove that the normaliser of any
finitely generated subgroup of a limit group is finitely generated.
Throughout this section, we work with rational homology. 
 However,
the proofs can easily be adapted to give analogous
 results for homology
with coefficients in any field, or in $\Z$.

\begin{lemma}\label{rk2} If $\G$ is a residually free group and 
 $S\subset\G$ is a non-cyclic subgroup, then
$H_1(S,\Q)$ has dimension at least 2 over $\Q$.
\end{lemma} 

\begin{proof} The lemma is obvious if  $S$ is abelian, because $\Gamma$
is torsion-free and contains no infinitely divisible elements.
If $S$ is non-abelian, then there exist $s,t\in S$ that do not commute. Since $\G$
is residually free, there exists a homomorphism
 $\phi:\G\to F$ with $F$ free and $\phi([s,t])\neq 1$.
Thus $S$ maps onto a non-abelian free group, namely
$\phi(S)$, and hence onto a free abelian group of rank
at least $2$.
\end{proof}

Recall the statement of Theorem \ref{1factor}.

\medskip

\noindent{\bf{Theorem 2}}
{\em{ Let $\Gamma$
  be a limit group and let $S\subset \Gamma$ be a subgroup. If
$H_1(S,\Q)$ has finite $\Q$-dimension,
then  $S$ is finitely generated (and hence a limit group).}}
\medskip

\begin{proof} Without loss 
of generality we may assume that $\G$ is
an $\o$-rft group, of height $h$ say. 
If $h=0$ then $S$ is a free product of free groups, 
surface groups and free abelian groups, and the result
is clear.  We may assume therefore that $h>1$. 
We consider the action of $S$ on the Bass-Serre tree
of the tower, pass to a {\em{minimal}} $S$-invariant subtree and
 form the quotient 
graph of groups. 
Thus we obtain a graph of groups decomposition $S=\pi_1(\mathcal{S},X)$
of $S$, with $X$ bipartite, in 
which the
edge groups are cyclic and the vertex groups are of two kinds
corresponding to the bipartite partition of vertices
(cf.~Lemma \ref{unconflate}): type (i)
vertex groups are the
intersections of $S$ with conjugates of 
a limit group $B\subset \G$
of height $(h-1)$;
type (ii) are the intersections of $S$ with the conjugates
of the free  or  free-abelian group at the top of the 
tower for $\G$. 

It is enough to prove the theorem for freely indecomposable
subgroups, so we may assume that all of the edge groups
are infinite cyclic. 

\smallskip

\noindent{\em Claim: If the graph $X$ is finite, then $S$ is finitely
generated.}

To see this, note that $H_1(S_e,\Q)\cong\Q$ for each of the finitely
many edges $e$ of $X$, so the Mayer-Vietoris sequence for $S$
shows that $H_1(S_v,\Q)$ is finite-dimensional for each vertex
$v$ of $X$.  By inductive hypothesis, each $S_v$ is finitely
generated, and hence so is $S$.

\smallskip
It remains to prove that $X$ is finite.  Certainly, since 
$H_1(S,\Q)$ is finite-dimensional, there is a finite connected subgraph
$Y$ of $X$, such that the inclusion-induced map
$$H_1(\pi_1(\mathcal{S},Y),\Q)\to H_1(\pi_1(\mathcal{S},X),\Q)=H_1(S,\Q)$$
is surjective.

In particular, $Y$ must contain all the simple  closed paths in $X$, so
all edges of $X\setminus Y$ separate $X$.  Suppose that $e$ is such an edge.  Let $Z_0$ be the component of $X\setminus\{e\}$ that contains
$Y$, and $Z_1$ the other component.  Then $e$ gives rise to a
splitting of $S$ as an amalgamated free product
$$S=S_0\ast_C S_1$$
with $C=S_e$ infinite cyclic, and $H_1(S_0,\Q)\to H_1(S,\Q)$ surjective.

The Mayer-Vietoris sequence of this decomposition shows that
$H_1(C,\Q)\to H_1(S_1,\Q)$ is surjective, so $H_1(S_1,\Q)$
has dimension at most $1$.  Since $\Gamma$ is
residually free,  Lemma \ref{rk2} tells us that $S_1$ must be
cyclic.  It follows that all the vertex groups of  $Z_1$ are cyclic.

By the argument given in the proof of Lemma \ref{BSsub}, the minimality
of $X$ means that vertex groups 
in $X$ of at most one type can be cyclic (type (i) if the top o
f the tower for $\Gamma$
is a quadratic block; type (ii) if it is an abelian block). Hence
$Z_1$ consists of a single vertex, $u$ say.
Moreover, if the top of the tower is an abelian block, then
$u$ is of type (ii), and $S_e$ is a maximal cyclic subgroup
of $S_u$, whence $S_e=S_u$, which contradicts the minimality of $X$.

\smallskip
We are reduced to the case where the top of the tower is a quadratic
block, and $X$ consists of the finite graph $Y$, together with
a (possibly infinite) collection of type (i) vertices of valence
$1$, each attached to a type (ii) vertex of $Y$ and each with
infinite cyclic vertex group. It follows easily that
$H_1(\pi_1(\mathcal{S},Y),\Q)\to H_1(S,\Q)$ is an isomorphism,
so $H_1(\pi_1(\mathcal{S},Y),\Q)$ is finite-dimensional.  Since
$Y$ is a finite graph and the edge groups are cyclic, it follows from
the Mayer-Vietoris sequence for $\pi_1(\mathcal{S},Y)$ that
$H_1(S_v,\Q)$ is finite-dimensional for each vertex $v$ of $Y$.

 To complete
the proof that $X$ is finite, we need only show that a type (ii)
vertex $v$ of $Y$ can have only finitely many neighbours in $X$.
To see this, note that $S_v$ is a surface group,
and the  edge-groups of the edges incident at $v$ form a subset
of the nontrivial peripheral subgroups of $S_v$. 
Thus $v$ is
incident to at most $1+\dim_\Q(H_1(S_v,\Q))<\infty$
edges of $X$, as required.
\end{proof}

\section{Normalisers are finitely generated}\label{nfg}

In this section we prove Theorem \ref{norm}:
%

\medskip

\noindent{\bf{Theorem 1}}
{\em{ If $\G$ is a limit group
and $H\subset\G$ is a finitely generated subgroup,
then either $H$ has finite index in its normaliser
or else the normaliser of $H$ is abelian.}}
\medskip

Since  subgroups of  fully residually free groups
are fully residually free,  finitely generated subgroups of
limit groups are limit groups. Thus
Theorem \ref{norm} is an immediate consequence
 of Theorem \ref{Lfi} and the following result.
 
\begin{thm}\label{fg} If $\G$ is a limit group
and $H\subset\G$ is a finitely generated subgroup,
then the normaliser in $G$ of
$H$ is finitely generated.
\end{thm}

\begin{proof} Let $X$ be a finite generating set for $H$
and let $N$ be the normaliser of $H$ in $\G$. All abelian
subgroups of limit groups are finitely generated \cite{BF},
so we assume that $N$ is not abelian.

Each finitely
generated subgroup $M\subset N$
is a  limit group, so if $X\subset M$, then
$H$ has finite index in $M$, by
Theorem \ref{Lfi}.
In particular, the dimension of $H_1(M,\Q)$
is bounded above by that of
 $H_1(H,\Q)$. Since $N$ is the union of its
finitely generated subgroups
 and homology commutes with direct
limits,  $H_1(N,\Q)$ has finite dimension. Hence $N$
is finitely generated, by
Theorem \ref{1factor}.
\end{proof}
 
\section{A splitting theorem}\label{spl}

The purpose of this section is to prove Theorem
\ref{split}.  In the statement we do not assume 
that $S$ is finitely generated.
%

\smallskip

\noindent{\bf{Theorem 3}}{\em{
Suppose that $\G_1,\dots,\G_n$ are limit groups 
and let $$S\subset \G_1\times\cdots\times \G_n$$ be an arbitrary subgroup.
If $L_i=\G_i\cap S$ is non-trivial and 
finitely generated for $i\le r$, then
a subgroup $S_0\subset S$ of finite index 
splits as a direct product
$$
S_0= L_1\times\dots\times L_r\times (S_0\cap G_r),
$$
where $G_r = \G_{r+1}\times\dots\times \G_n$.}}

\bigskip

\begin{proof} Let $P_i$ be the projection of $S$ to $\G_i$. 
Since $L_i$ is normal in $S$, it is normal in $P_i$. Finitely
generated subgroups of limit groups are limit groups, so it follows
from  Theorem \ref{Lfi} that $L_i$ has finite index in $P_i$.

We have shown that $D=L_1\times\dots\times L_r$ has finite index in the
projection of $S$ to $\G_1\times\dots \G_r$. Thus the inverse
image $S_0\subset S$ of $D$ is of finite index, 
and by splitting the short
exact sequence $1\to (S_0\cap \G_r)\to S_0\to D\to 1$ of the
projection we obtain the desired direct product decomposition of $S_0$.
\end{proof}

\medskip\noindent{\bf Acknowledgement:}  
We thank Mladen
Bestvina, Zlil Sela and Henry Wilton
for comments concerning this work.

\medskip\centerline{\bf Authors' addresses}

\smallskip\begin{center}\begin{tabular}{lll}
Martin R. Bridson &\qquad\qquad & James Howie\\
Department of Mathematics && Department of Mathematics\\
Imperial College London && Heriot-Watt University\\
London SW7 2AZ && Edinburgh EH14 4AS\\
{\tt m.bridson@imperial.ac.uk} && {\tt J.Howie@hw.ac.uk}
\end{tabular}\end{center}

\end{document}